\documentclass[a4paper,twoside,12pt]{article}

\newtheorem{thm}{Theorem}

\usepackage{amsmath} 
\usepackage[english]{babel}
\usepackage{babelbib}
\usepackage[authoryear,square,comma]{natbib}
\RequirePackage{geometry}
\usepackage[utf8]{inputenc}
\usepackage[T1]{fontenc}
\setcitestyle{square,aysep={}}

\usepackage{comment}
\usepackage{float}
\usepackage{etoolbox}

\usepackage{txfonts} 

\usepackage[pdftex]{graphicx}

\usepackage{listingsutf8}

\title{On the extension of $D(-8k^2)$-pair $\{8k^2, 8k^2+1\}$ }
\author{Nikola Adžaga\footnote{Faculty of Civil Engineering, University of Zagreb, e-mail: nadzaga@grad.hr}, Alan Filipin\footnote{Faculty of Civil Engineering, University of Zagreb, e-mail: filipin@grad.hr}}
\date{} 

\usepackage{xcolor}

\newcounter{grammarchild}
\usepackage{hyperref}
\hypersetup{
 unicode=true,
 colorlinks=true,
 linkcolor=black,
 citecolor=black,
 filecolor=black,
 urlcolor=black,
 bookmarksopen=true,
 bookmarksopenlevel=0,
 bookmarksnumbered=true,
 pdftitle={Automated conjecturing of Frobenius numbers via grammatical evolution},
 pdfauthor={Nikola Adžaga},
 pdfkeywords={},
 pdfstartpage={},
 pdfstartview=FitH,
 pdfnewwindow=true
}

\begin{document}
\maketitle

\section*{Abstract}
Let $n$ be a nonzero integer. A set of $m$ positive integers is called a $D(n)$-$m$-tuple if the product of any two of its distinct elements increased by $n$ is a perfect square. Let $k$ be a positive integer. By elementary means, we show that the $D(-8k^2)$-pair $\{8k^2, 8k^2+1\}$ can be extended to at most a quadruple (the third and fourth element can only be $1$ and $32k^2+1$). At the end, we suggest considering a $D(-k^2)$-triple $\{ 1, 2k^2, 2k^2+2k+1\}$ as possible future research direction.


\section{Introduction}
The research of Diophantine m-tuples has a long history. From Diophantus, through Fermat and Euler, this field has remained interesting to number theorists. In 20th century, Baker developed theory of linear forms in logarithms which he used to prove the first results on nonextendibility of such sets (\cite{baker}).

A set of $m$ positive integers $\{a_1, a_2, \dots, a_m\}$ is called $D(n)$-$m$-tuple if $a_i a_j+n$ is a perfect square for all $1 \leqslant i < j \leqslant m$. Natural question regarding this sets is their possible size. \textcolor{black}{Most of the results in this field, especially newer ones, rely on linear forms in logarithms (e.g.~\cite{BFT}), tools from Diophantine approximation (such as \cite{Bennett}) and elliptic curves (\cite{sestorke}).  }

On the other hand, Fujita \& Togbe (\cite{Fujita2012}) had proven\textcolor{black}{, in an elementary, and relatively simple manner,} if $\{k^2, k^2+1, c, d\}$ is a $D(-k^2)$-quadruple with $c < d$, then $c=1$ and $d=4k^2+1$ (in that case, $3k^2+1$ must be a square).

Similarly, $\{1, 2k^2, 2k^2+1\}$ is a $D(-2k^2)$-triple. However, for odd $k$, the number added, $-2k^2$, leaves a remainder of $2$ when divided by $4$, and already in 1985, Brown (\cite{brown}) has shown that there is no $D(n)$-quadruple for $n \equiv 2 \pmod{4}$, simply by observing quadratic residues modulo $4$.


Therefore, this problem is interesting only for even $k$. We rephrase it in the following manner. Let $k$ be a positive integer. Set $\{8k^2, 8k^2+1\}$ has $D(-8k^2)$-property, as well as the set $\{1, 8k^2, 8k^2+1\}$, i.e.~the product of any two of its distinct elements subtracted by $8k^2$ is a perfect square:
\[ 1\cdot 8k^2 - 8k^2 = 0, \quad 1\cdot (8k^2+1)-8k^2 = 1,\quad  8k^2 \cdot (8k^2+1) - 8k^2 = 64k^4. \]

The main result is the following theorem.
\begin{thm}
\label{par}
Let $k$ be a positive integer. If a set $\{8k^2, 8k^2+1, c, d\}$ is a $D(-8k^2)$-quadruple such that $d > c$, then $c=1$ and $d=32k^2+1$. In this case, $k$ can be expressed as $k = k_i$, where $k_i$ is the recurrent sequence defined by $k_1 = 1, k_2=10, k_{n+2} = 10k_{n+1}-k_n$. 
\end{thm}

We prove this by first considering the extendibility of the triple $\{1, 8k^2, 8k^2+1\}$.

\section{The extendibility of the $D(-8k^2)$-triple $\{1, 8k^2, 8k^2+1\}$}
\label{sec:triple}
In this section, we prove the following theorem
\begin{thm}
\label{trojka}
If a set $\{1, 8k^2, 8k^2+1, d\}$ is a $D(-8k^2)$-quadruple, then $d=32k^2+1$. In that case, $24k^2+1$ must be a square.
\end{thm}

Extending the initial triple with $d$ and then eliminating $d$ leads to a system consisting of a Pell ($z^2-(16k^2+2)y^2=1$) and a pellian equation ($x^2-2y^2=-8k^2+1$). By solving Pell equation, we get two recurrent sequences $y_n$ and $z_n$. Due to the second equation, the problem reduces to examining when can an element of the new sequence
$X_n = 2y_n^2-8k^2+1$ be a complete square. Using the relations between $y_n$ and $z_n$, e.g.~$y_{2n+1} = 2y_nz_n$, we write $X_n$ as a product of two factors, one of which is obviously not a square. We finish the proof by showing that these factors are relatively prime via principle of descent.

\subsection{System of simultaneous pellian equations}
\label{subsec:system}
If the set $\{ 1, 8k^2, 8k^2+1, d \}$ has $D(-8k^2)$-property, then there exist integers $x, y', z$ such that
\begin{align*}
d-8k^2=x^2, \\
8k^2d-8k^2 = (y')^2,\\
(8k^2+1)d-8k^2 = z^2
\end{align*}
In the second equation we note that $8k^2$ divides $(y')^2$, hence $4k$ divides $y'$, i.e.~$y'=4ky$ for some integer $y$. Dividing this equation with $8k^2$, we obtain a simpler one, $d-1 = 2y^2$. From this and from the first equation we easily eliminate $d$ to obtain a pellian equation $x^2-2y^2 = -8k^2+1$. In a similar manner, by eliminating $d$ from the third equation and the simplified second equation, we obtain Pell's equation $z^2-(16k^2+2)y^2 = 1$.

To summarize, if a $D(-8k^2)$-triple can be extended, then there exist integers $x, y$ and $z$ satisfying the following system:
\begin{align}
x^2-2y^2 = -8k^2+1 \label{1} \\
z^2-(16k^2+2)y^2 = 1 \label{2}
\end{align}

\textcolor{black}{This system is more complicated than the one obtained in \cite{Fujita2012}, as one of the equations there was very simple ($y^2-x^2=k^2-1$).}

\textcolor{black}{The equation \eqref{1} is a pellian equation, and it can have a large number of infinite classes of solutions. On the other hand, }Pell's equation (\ref{2}) has one (infinite) class of solutions, generated by a fundamental one. Using continued fractions ($\sqrt{16k^2+2} = [4k; \overline{4k, 8k}]$), we obtain its fundamental solution $16k^2+1+4k\sqrt{16k^2+2}$. Thus we easily obtain recurrent system of sequences $(z_n)_{n \geqslant -1}$ and $(y_n)_{n \geqslant -1}$ containing all nonnegative integer solutions $(z, y)$ to (\ref{2}):
\begin{align}
 z_{n+1} &= (16k^2+1)z_n + 4k(16k^2+2)y_n, \quad z_0 = 16k^2+1,\, z_{-1} = 1 \label{rec1} \\
 y_{n+1} &=(16k^2+1)y_n + 4k z_n, \quad \quad \quad \quad \quad \,\,\,y_0 = 4k, \,y_{-1} = 0. \label{rec2}
\end{align}

From this, one can get second order linear recurrence relations for $(y_n)$ and $(z_n)$, solve them and obtain explicit expressions for all solutions of (\ref{2}). This, along with the proofs of some identities relating sequences $(y_n)$ and $(z_n)$, is done in the appendix.

On the other hand, $y$ should also be a solution of (\ref{1}), i.e.~$x^2$ should be $2y_n^2-8k^2+1$ for some $n \in \mathbb{N}_0$. This is why a new sequence is introduced with the following definition $X_n := 2y_n^2-8k^2+1$.

\subsection{Squares in a sequence $X_n$}
 Using the relations between $y_n$ and $z_n$, e.g.~$y_{2n+1} = 2y_nz_n$, we write $X_n$ as a product of two factors, one of which is obviously not a square. We show that these factors are relatively prime (via principle of descent for odd indices). However, the factorization and the whole proof depends on the oddness of the index.

\subsubsection{Odd indices}
For odd indices, we use the following identity:
\begin{align}
 y_{2n+1} = 2y_n z_n, \label{idodd}
\end{align}
to show that $X_{2n+1}$ is never a square. This identity (\ref{idodd}) is proven in the appendix. It enables us to do the following computation:

\begin{align*}
X_{2n+1} &= 2y_ {2n+1}^2-8k^2+1 \\
&= 8y_n^2z_n^2-8k^2+1  \,\, (\text{substitute }z_n^2 \text{ from (\ref{2})})\\
&= 8y_n^2 (1+(16k^2+2)y_n^2) - 8k^2+1\\
&= 8y_n^2+16(8k^2+1)y_n^4 - 8k^2 + 1 \,\, (\text{can be factored as})\\
&=  (4y_n^2+1)(32y_n^2 k^2 + 4y_n^2 - 8k^2 + 1).
\end{align*}
The first factor cannot be a square, since $(2y_n)^2 < 4y_n^2+1 < (2y_n+1)^2$ (for $y_n \in \mathbb{N}$). Therefore, it suffices to show that the two factors obtained are relatively prime. We prove this via principle of descent.

More precisely, we prove that
\begin{align*}
 p \mid 4y_n^2+1 & \text{ and } p\mid 32y_n^2 k^2 + 4y_n^2 - 8k^2 + 1  \Longrightarrow \\
\Longrightarrow p \mid 4y_{n-1}^2+1 & \text{ and } p\mid 32y_{n-1}^2 k^2 + 4y_{n-1}^2 - 8k^2 + 1, \\
\end{align*}
\vskip -2em
which will lead us into a contradiction.

Assume that prime $p$ divides both $4y_n^2+1$ and $32y_n^2 k^2 + 4y_n^2 - 8k^2 + 1$. Then $p$ divides $32y_n^2 k^2 + 4y_n^2 - 8k^2 + 1 - (4y_n^2+1) = 8k^2(4y_n^2 - 1)$ too. Since $p$ is odd (because it divides an odd number $4y_n^2+1$), it follows that $p \mid k^2(4y_n^2 - 1)$. It cannot divide a second factor, because then it would divide $4y_n^2+1 - (4y_n^2 - 1) = 2$ as well. To conclude, $p$ divides $k^2$, but then it also divides $k$.

Let us prove now that $p$ divides both $4y_{n-1}^2+1$ and $32y_{n-1}^2 k^2 + 4y_{n-1}^2 - 8k^2 + 1$. From (\ref{rec2}), recurrence relation for $(y_n)$, we get $y_n - y_{n-1} = 16k^2 y_{n-1}+4kz_{n-1}$, and since $p$ divides $k$, it follows that it divides the right hand side, which is equal to $y_n - y_{n-1}$. Hence, $y_n  \equiv y_{n-1} \pmod{p}$. Indeed, $p$ divides $4y_{n-1}^2+1$. On the other hand, since $32y_{n-1}^2 k^2 + 4y_{n-1}^2 - 8k^2 + 1 = 8k^2(4y_{n-1}^2 - 1) + 4y_{n-1}^2+1$, it is true that $p$ divides $32y_{n-1}^2 k^2 + 4y_{n-1}^2 - 8k^2 + 1$.

Further descent implies that $p$ divides $4y_0^2 + 1 = 64k^2+1$ (and $32y_0^2 k^2 + 4y_0^2 - 8k^2 + 1$). However, since $p$ divides $k$, it would divide $1$,  which is a contradiction. To conclude, $4y_n^2+1$ and $32y_n^2 k^2 + 4y_n^2 - 8k^2 + 1$ do not have any common prime factors, i.e.~they are relatively prime.

\subsubsection{Even positive indices}
For even indices, we use the following identity:
\begin{align}
z_{2n}-1 = \left( \frac{y_n + y_{n-1}}{4k} \right)^2, \label{ideven}
\end{align}
to show that $X_{2n}$ is not a square for positive integer $n$. We easily compute
\begin{align*}
X_{2n} &=2y_{2n}^2-8k^2 + 1  \stackrel{(\ref{2})}{=} 2\cdot \frac{z_{2n}^2-1}{16k^2+2} - 8k^2 + 1= \\
	&= \frac{z_{2n}^2-1-64k^4-8k^2+8k^2+1}{8k^2+1} \\ 
	&= \frac{z_{2n}^2-64k^2}{8k^2+1} \\
	&= \frac{z_{2n}-8k^2}{8k^2+1} \cdot (z_{2n}+8k^2),
\end{align*}
We show that the second factor, $z_{2n}+8k^2$ cannot be a square for $n\geqslant 1$. Due to the identity \eqref{ideven}, we obtain \[\left( \frac{y_n + y_{n-1}}{4k} \right)^2 < z_{2n}+8k^2 =  \left( \frac{y_n + y_{n-1}}{4k} \right)^2 +1 + 8k^2 < \left(\frac{y_n + y_{n-1}}{4k} + 1\right)^2. \]

Last inequality holds if and only if $\displaystyle 1+8k^2 < \frac{y_n + y_{n-1}}{2k} + 1$, i.e.~ if and only if $16k^3 < y_n + y_{n-1}$.

Since already $y_1 = 2(16k^2+1) \cdot 4k$ \textcolor{black}{and $y_n$ is obviously an increasing sequence, the required} inequality truly holds for $n  \geqslant 1$.

In this case, it is easier to prove that the factors obtained are relatively prime. Assume prime $p$ divides both $\displaystyle \frac{z_{2n}-8k^2}{8k^2+1}$ and $z_{2n}+8k^2$. Then it divides $z_{2n}-8k^2$, as well as the difference $z_{2n}+8k^2 - (z_{2n}-8k^2) = 16k^2$. Since $p$ divides $z_{2n}+8k^2$, and all elements of a sequence $z_n$ are odd (see (\ref{rec1})), it follows that $p$ is odd. Hence $p\mid 16k^2$ implies $p\mid k^2$. It follows that $p\mid z_{2n}$. But the sequence $z_n$ is relatively prime with $k$ since all elements give a remainder of $1$ when divided with $k$. Therefore, it is impossible that $z_{2n}$ and $k$ have a common prime factor $p$, so we obtain a contradiction again.

One should also check that the first factor, fraction $\displaystyle \frac{z_{2n}-8k^2}{8k^2+1}$, is actually an integer. But that is easy: inductively prove that elements of $z_n$, when divided by $8k^2+1$, leave remainders $( 1, 8k^2, 1, 8k^2, \dots)$. Hence $z_{2n}-8k^2$ is indeed divisible by $8k^2+1$.

We are left with the case $n = 0$, i.e.~when $\displaystyle X_0 = \frac{z_{0}-8k^2}{8k^2+1} \cdot (z_{0}+8k^2) = 24k^2+1$ is a square.

\subsubsection{Conclusion for the extendibility of a triple}
To conclude, the only element of a sequence $(X_n)_{n \geqslant 0}$ which can be a perfect square is $X_0 =  24k^2+1$. In that case, since $x^2=24k^2+1 \label{P}$, we see that the fourth element can only be $d=32k^2+1$.  However, $24k^2+1$ is not a square for all $k$, \textcolor{black}{but to determine for which $k$ it is a perfect square, we just have to solve another Pell's equation $m^2-24k^2=1$. }

\section{The extendibility of the $D(-8k^2)$-pair $\{8k^2, 8k^2+1\}$}
In this section, we prove theorem \ref{par}.

Assume that $\{8k^2, 8k^2+1, c\}$ is a $D(-8k^2)$-triple, where $1 < c$. Then $8k^2 c - 8k^2 = (s')^2$, hence $s' = 4ks$, which simplifies the equation to $c-1=2s^2$. Also, $(8k^2+1)c-8k^2 = t^2$ ($t$ and $s$ are nonnegative integers). Eliminating $c$, we obtain $t^2-(16k^2+2)s^2 = 1$. From \ref{subsec:system}, we already know that $t+s\sqrt{16k^2+2} = (16k^2+1+4k\sqrt{16k^2+2})^\nu$ for $\nu \in \mathbb{N}_0$. Hence, $s = s_\nu$ for some $\nu$, where $s_0 = 0, s_1 = 4k, s_{\nu+2} = 2(16k^2+1)s_{\nu+1}-s_\nu$.

Assume now that there exist $D(-8k^2)$-quadruples $\{8k^2, 8k^2+1, c, d\}$ such that $1<c, d$, i.e.~let integers $d > c > 1$ be such that $\{8k^2, 8k^2+1, c, d\}$ is a $D(-8k^2)$-quadruple with minimal $c > 1$. Therefore \[8k^2d-8k^2 = (x')^2, \quad (8k^2+1)d-8k^2 = y^2, \quad cd - 8k^2 = z^2.\] First equation simplifies again to $d-1 = 2x^2$. Eliminating $d$, we obtain a system
\begin{align}
y^2-(16k^2+2)x^2 &= 1 \\
z^2-2cx^2 &= c-8k^2 \label{2.2}
\end{align}

and again $y+x\sqrt{16k^2+2} = (16k^2+1+4k\sqrt{16k^2+2})^m, \quad m \in \mathbb{N}_0$, i.e.\ $x$ is an element of the sequence $v_0 = 0, v_1 = 4k, v_{m+2} = 2(16k^2+1)v_{m+1}-v_m$.

On the other hand, fundamental solution of Pell's equation ($z^2-2cx^2=1$), corresponding to pellian equation \eqref{2.2}, is $(2c-1, 2s)$. (since $(2c-1)^2-2c(2s)^2 = 4c^2-4c+1-4c(c-1) = 1$). 
Therefore, if $(z, x)$ is a solution of \eqref{2.2}, then there exists a fundamental solution $(z_0, x_0)$ such that $z+x\sqrt{2c} = (z_0+x_0\sqrt{2c})(2c-1+2s\sqrt{2c})^n$, and \[ 0 < x_0 \leqslant \frac{2s}{\sqrt{2(2c-2)}} \sqrt{c-8k^2} = \frac{s}{\sqrt{2s^2}} \sqrt{c-8k^2} = \sqrt{\frac{c}{2} - 4k^2} = \sqrt{s^2+\frac 12 - 4k^2} < s .\]

(The second inequality above comes from Theorem $108\text{a}$ in \cite{nagell}).

We conclude that $x$ is also an element of a sequence $w_0 = x_0, w_1 = (2c-1)x_0 + 2sz_0, w_{n+2} = 2(2c-1)w_{n+1}-w_n$. Observe that $w_n \equiv x_0 \pmod{s}$ for all $n$, since $c = 2s^2+1 \equiv 1 \pmod{s}$. Given that
\[ (v_m)_{m \geqslant 0} \equiv (0, s_1, s_2, \dots, s_{\nu-1}, 0, -s_{\nu-1}, -s_{\nu-2}, \dots, -s_1, 0, s_1, \dots) \pmod{s},\] and $x_0 < s$, it follows that $x_0^2 = s_i^2$ for some $i < \nu$.

Let now $d_0 = 2x_0^2+1$. Then $(8k^2+1)d_0-8k^2 = (16k^2+2)x_0^2 + 1 = (16k^2+2)s_i ^2 + 1$, so \[ (8k^2+1)d_0-8k^2 = t_i^2. \] Also, \[ 8k^2d_0-8k^2 = 16k^2x_0^2 = (4kx_0)^2 \] and \[ cd_0 - 8k^2 = c(2x_0^2+1)-8k^2 = 2cx_0^2+c-8k^2 = z_0^2-c+8k^2+c-8k^2 = z_0^2. \]

Therefore, a set $\{8k^2, 8k^2+1, d_0, c\}$ is also a $D(-8k^2)$-quadruple. Observe that $d_0 = 2x_0^2+1 < 2s^2+1 = c$, and, since $c$ is the minimum possible value greater than $1$, $d_0$ must be $d_0=1$, i.e.~$x_0 = 0 > 0$, which is a contradiction.

There is no $D(-8k^2)$-quadruple $\{8k^2, 8k^2+1, c, d\}$ such that $1 < c < d$.

For the last assertion, we use theorem \ref{trojka}, to conclude that $24k^2+1$ must be a square. But this is again Pell's equation: $m^2-24k^2=1$. Fundamental solution $5+\sqrt{24}$ gives us the recurrent sequence for $k$, $k_1= 1, k_2 = 10, k_{n+2} = 10k_{n+1}-k_n, \, \forall n \in\mathbb{N}$.

\section{Similar problem}
We can consider $D(-k^2)$-triple $\{ 1, 2k^2, 2k^2+2k+1\}$ in a similar manner. Extending it with $d$, we obtain the following system of simultaneous pellian equations:
\begin{align*}
y^2-2x^2 = 2k^2-1 \\
z^2-(4k^2+4k+2)y^2 = 4k+2
\end{align*}
Even though none of this equations is Pell's, second equation has only one fundamental solution, $(z^*, y^*)=(2k+2, 1)$. We show this by placing $z^*$ between squares of $(2k+1)y$ and $(2k+1)y+1$, using the upper bound for $y^*$ (found in \cite{nagell} again).

This gives us recurrences for the solutions of the second equation, $y_{n+1} = (8k^2+8k+3)y_n+(4k+2)z_n$ from which we conclude that $y_{n+1} \equiv y_n \pmod{4k+2}$, i.~e.~$y_n \equiv y_0=1 \pmod{4k+2}$. Similarly, if we change one of the signs in fundamental solution, say \\$(2k+2, -1)$, to generate solutions $(z_n', y_n')$, we can conclude that $y_n' \equiv y_0'= -1 \pmod{4k+2}$.
In any case, from the first equation, $2x^2=y^2-2k^2+1 \equiv 1-2k^2+1=2(1-k^2) \pmod{4k+2}$, so $x^2 \equiv 1-k^2 \pmod{2k+1}$ and \[ x^2 \equiv 1-k^2 = (1-k)(1+k) \equiv (k+2)(k+1)=k^2+3k+2\equiv k^2+k+1 \pmod{2k+1} \]
Multiplying by $4$, $4x^2 \equiv 4k^2+4k+4-(2k+1)^2 = 3 \pmod{2k+1}.$

Hence, $3$ is a square modulo $2k+1$. It is a quadratic residue modulo every prime factor of $2k+1$ as well. Since $3$ is a quadratic residue modulo prime $p$ if and only if $p \equiv \pm 1 \pmod{12}$, it follows that all prime factors of $2k+1$ are of that form, or equal to $3$. Only one $3$ can divide it, since $3\mid 4x^2-3$ implies $9\nmid 4x^2-3$.

We can conclude that for $k\equiv 2, 3 \pmod{6}$, there is no extension. However, we conjecture that the fourth element can only be $d=8k^2+4k+1$, and in that case, $7k^2+4k+1$ must be a perfect square.

\section*{Acknowledgments}
The authors are supported by Croatian Science Foundation under the project no. 6422.

\section*{Appendix}
In Section \ref{sec:triple}, we used the fact that, if we denote by $(z_n, y_n)_{n\geqslant -1}$ all nonegative integral solutions to $z^2-(16k^2+2)y^2=1$, arranged in ascending order, then $y_{2n+1} = 2y_n z_n$. Here we provide the proof of this identity. Fundamental solution of this Pell's equation is $16k^2+1+ 4k\sqrt{16k^2+2}$, which means that every solution can be obtained from the previous one by $z_{n+1}+y_{n+1}\sqrt{16k^2+2} = (z_n+y_n\sqrt{16k^2+2})(16k^2+1+ 4k\sqrt{16k^2+2})$. This gives us mixed recurrences (\ref{rec1}) and (\ref{rec2}). Now
 \begin{align*}
z_{n+2} &= (16k^2+1)z_{n+1}+4k(16k^2+2)\cdot((16k^2+1)y_n+4kz_n)\\
	&= (16k^2+1)z_{n+1}+16k^2(16k^2+2)z_n + (16k^2+1)(z_{n+1}-(16k^2+1)z_n)\\
	&= 2(16k^2+1)z_{n+1} - z_n,
\end{align*}
with initial conditions $z_{-1} = 1, z_0 = 16k^2+1$.

In the same manner, we get second order linear recurrence relation for $(y_n)$: $y_{n+2} = 2(16k^2+1)y_{n+1} - y_n$, with initial conditions $y_{-1} = 0, y_0 = 4k$.
Solving these recurrences, we get explicit expressions:\begin{align*}
 y_n &=c_1 (16k^2+1+4k\sqrt{16k^2+2})^n + c_2 (16k^2+1-4k\sqrt{16k^2+2})^n \\
z_n &= c_3 (16k^2+1+4k\sqrt{16k^2+2})^n + c_4 (16k^2+1-4k\sqrt{16k^2+2})^n,
\end{align*} where \[ c_1 = \frac{64k^3+8k+(16k^2+1)\sqrt{16k^2+2}}{2(16k^2+2)}, \quad c_2 = \frac{64k^3+8k-(16k^2+1)\sqrt{16k^2+2}}{2(16k^2+2)} \]
\[ c_3 = \frac{128 k^4 + 24 k^2 + 1 + (32 k^3 + 4 k)\sqrt{16k^2+2}}{16k^2+2}, \quad c_4 = \frac{128 k^4 + 24 k^2 + 1 - (32 k^3 + 4 k)\sqrt{16k^2+2}}{16k^2+2} \]

We compute
\begin{align*}
2y_n z_n &= 2(c_1 (16k^2+1+4k\sqrt{16k^2+2})^n + c_2 (16k^2+1-4k\sqrt{16k^2+2})^n) \cdot \\
	& \quad \cdot (c_3 (16k^2+1+4k\sqrt{16k^2+2})^n + c_4 (16k^2+1-4k\sqrt{16k^2+2})^n) = \\
&= 2\left(c_1c_3(16k^2+1+4k\sqrt{16k^2+2})^{2n} + c_2c_4(16k^2+1-4k\sqrt{16k^2+2})^{2n}+c_1c_4+c_2c_3\right)
\end{align*}

Now all we need to do is compute, i.e.~compare these coefficients (this is best done by computer):
 \begin{align*}
2c_1c_3 &= 2 \frac{64k^3+8k+(16k^2+1)\sqrt{16k^2+2}}{2(16k^2+2)} \cdot\frac{128 k^4 + 24 k^2 + 1 + (32 k^3 + 4 k)\sqrt{16k^2+2}}{16k^2+2} = \\
&= \frac{(64 k^3 + 8 k) (128 k^4 + 24 k^2 + 1) + (16 k^2 + 1) (32 k^3 + 
    4 k) (16 k^2 + 2)  } {(16k^2+2)^2} +\\
&\quad + \frac{((64 k^3 + 8 k) (32 k^3 + 4 k) + (16 k^2 + 1) (128 k^4 + 24 k^2 + 1))\sqrt{16k^2+2}}{(16k^2+2)^2} = \\
&= \frac{16 k + 512 k^3 + 5120 k^5 + 16384 k^7 + (1 + 72 k^2 + 1024 k^4 + 4096 k^6)\sqrt{16k^2+2}}{(16k^2+2)^2} =\\
&= \frac{16 k (8 k^2+1)^2 (16 k^2+1) + (8k^2+1)(512k^4+64k^2+1)\sqrt{16k^2+2}}{(16k^2+2)^2} =\\
&= \frac{8k(8k^2+1)(16k^2+1)+\frac 12 (512k^4+64k^2+1)\sqrt{16k^2+2}}{16k^2+2} = \\
&= \frac{16k(8k^2+1)(16k^2+1)+(512k^4+64k^2+1)\sqrt{16k^2+2}}{2(16k^2+2)}
\end{align*}
On the other hand, in $y_{2n+1}$, coefficient of $(16k^2+1+4k\sqrt{16k^2+2})^{2n}$ is
\[
c_1 (16k^2+1+4k\sqrt{16k^2+2}) = \frac{64k^3+8k+(16k^2+1)\sqrt{16k^2+2}}{2(16k^2+2)}\cdot (16k^2+1+4k\sqrt{16k^2+2}) =\]
\[ = \frac{(64 k^3 + 8k)(16 k^2 + 1) + 4k(16k^2+1)(16k^2+2) + ((64 k^3 + 8k)\cdot 4k + (16k^2+1)^2)\sqrt{16k^2+2}}{2(16k^2+2)} = 
\]
\[ = \frac{16k(128k^4+24k^2+1) + (512k^4+64k^2+1)\sqrt{16k^2+2}}{2(16k^2+2)}, \]

so indeed the coefficient of $(16k^2+1+4k\sqrt{16k^2+2})^{2n}$ is equal in $y_{2n+1}$ and in $y_nz_n$.

Analogously for the conjugate $(16k^2+1-4k\sqrt{16k^2+2})^{2n}$, while $c_1c_4+c_2c_3 = 0$.

The other identity, $\displaystyle z_{2n}-1 = \left( \frac{y_n + y_{n-1}}{4k} \right)^2$, is proven in a similar manner.

\bibliographystyle{plainnat} 
\bibliography{brefs}

\end{document}